\documentclass[12pt]{amsart}

\usepackage{amsfonts, amscd}
\usepackage[dvips]{graphicx}
\author{Sergey Maksymenko}
\title[Circle-valued Morse functions]
{Connected components of the space of circle-valued Morse functions on surfaces}

\swapnumbers
\newcommand\globcnt{subsubsection}

\theoremstyle{plain}

\newtheorem{theorem}[\globcnt]{Theorem}
\newtheorem{lemma}[\globcnt]{Lemma}
\newtheorem{prop}[\globcnt]{Proposition}

\newtheorem{claim}[\globcnt]{Claim}
\theoremstyle{definition}

\newtheorem{defn}[\globcnt]{Definition}

\newcommand{\RRR}{{\mathbb R}}
\newcommand{\ZZZ}{{\mathbb Z}}

\newcommand{\Int}{\rm{Int}}
\newcommand{\grad}{\nabla} 

\newcommand{\eps}{\varepsilon}

\newcommand\MM{{\mathcal M}}

\newcommand\surf{M}

\newcommand\onemanif{P}
\newcommand\Circle{S^1}
\newcommand\Rline{\RRR^1}

\newcommand\MrsMP{\MM(\surf,\onemanif)}
\newcommand\MrsMR{\MM(\surf,\Rline)}
\newcommand\MrsMS{\MM(\surf,\Circle)}

\newcommand\smsp[1]{C^{\infty}\ifx#1\empty\else(#1)\fi}
\newcommand\smMP{\smsp{\surf,\onemanif}}

\newcommand\func{f}
\newcommand\gfunc{g}
\newcommand\hfunc{h}

\newcommand\Func{F}
\newcommand\cfunc{\tfunc}

\newcommand\sgmhom{$\Sigma$}

\newcommand\sgmh{\stackrel{\Sigma}{\sim}}

\newcommand\Interv{I}
\newcommand\aval{0}
\newcommand\bval{1}
\newcommand\imInterv{[\aval,\bval]}

\newcommand\comp{V}
\newcommand\compA{\comp_{\aval}}
\newcommand\compB{\comp_{\bval}}

\newcommand\crpnt[2]{c_{#1}(#2)}
\newcommand\crpntf[1]{\crpnt{#1}{\func}}

\newcommand\bcomp{V}
\newcommand\sgn[1]{\eps_{#1}}

\newcommand\Crtype[1]{K(#1)}

\newcommand\pnt{z}
\newcommand\apnt{x}
\newcommand\bpnt{y}
\newcommand\xpnt{s}

\newcommand\tfunc{\widetilde{\func}}
\newcommand\tgfunc{\widetilde{\gfunc}}
\newcommand\thfunc{\widetilde{\hfunc}}

\newcommand\pr{p}
\newcommand\qpr{q}

\newcommand\regcomp{X}
\newcommand\tregcomp{\widetilde{\regcomp}}
\newcommand\gregcomp{Y}
\newcommand\tgregcomp{\widetilde{\gregcomp}}
\newcommand\concomp{C}
\newcommand\aconcomp{D}
\newcommand\aconcomps{V}

\newcommand\tsurf{\widetilde{\surf}}

\newcommand\csurf[1]{ \widetilde{\surf} } 
\newcommand\ctsurf[1]{ \widehat{\surf} } 

\newcommand\cpr{\pr}

\newcommand\orpr{\tau}

\newcommand\orfunc{\widehat{\func}}
\newcommand\orgfunc{\widehat{\gfunc}}

\newcommand\orsurf{\widehat{\surf}}
\newcommand\orconcomp{\widehat{\concomp}}
\newcommand\oraconcomp{\widehat{\aconcomp}}

\newcommand\orregcomp{\widehat{\regcomp}}
\newcommand\orgregcomp{\widehat{\gregcomp}}

\newcommand\nbh{U}

\newcommand\crv{\gamma}
\newcommand\tcrv{\widetilde{\crv}}

\newcommand\larcD{L_{0}}
\newcommand\larcU{L^{1}}
\newcommand\larcDU{L_{0}^{1}}
\newcommand\larcC{L_{c}}

\begin{document}

\begin{abstract}
We classify the path-components of the space of circle-valued Morse functions on compact surfaces: two Morse functions $\func, \gfunc: \surf \to \Circle$ belong to same path-component of this space if and only if they are homotopic and have equal numbers of critical points at each index.
\end{abstract}

\maketitle


 \section{Introduction}
 \label{sect:Intro}
Let $\surf$ be a smooth ($C^{\infty}$) connected compact surface, orientable or not, with boundary $\partial\surf$ or without it, and $\onemanif$ a one-dimensional manifold either the real line $\Rline$ or the circle $\Circle$.
Let $\MrsMP$ denotes the subspace of $\smMP$ consisting of Morse mappings $\surf\to\onemanif$.
It is well-known (e.g. Milnor~\cite{Milnor:HCobTheory}) that for the case $\partial\surf=\emptyset$ the set $\MrsMP$ is an everywhere dense open in $\smMP$ with the $C^{\infty}$ Whitney topology of $\smMP$.

Recently, S.~V.~Matveev (his proof is included and generalized in
E.~Kud\-ryavt\-seva~\cite{Kudr}) and V.~V.~Sharko~\cite{Sharko}
have obtained a full description of the set $\pi_0\MrsMR$ of connected path-components of $\MrsMR$. 
Their methods are independent and based on different ideas.

For orientable closed surfaces the classification of $\pi_0\MrsMS$ was initially given in author's Ph.D, see~\cite{Maks}.
The problem was proposed to the author by V.~V.~Sharko.
In this note, we extend the results of~\cite{Maks} to all compact surfaces (Theorem~\ref{th:Main}) and simplify their proof.


To begin with, let us fix, once and for all, some Riemannian metric on $\surf$ and some orientation of $\onemanif$.

A $C^{\infty}$-mapping $\func:\surf\to\onemanif$ is {\em Morse\/} if the following conditions hold true:
\begin{enumerate}
\item
all critical points of $\func$ are non-degenerated and belong to the interior of $\surf$;
\item
$\func$ is constant at every connected component of $\partial\surf$
though its values on different components may differ each from other.
\end{enumerate}

Let $\func:\surf\to\onemanif$ be a Morse function and $\pnt$ a non-degenerated critical point $\func$. 
Define the \emph{index} of $\pnt$ to be the usual one with respect to an arbitrary local representation \ $\surf \supset \RRR^2 \stackrel{\func}{\to} \RRR^1 \stackrel{\phi}{\subset} \onemanif$ \ of $\func$ in which $\phi$ preserves orientations.
Denote by $\crpnt{i}{\func}$ $(i=0,1,2)$ the number of critical points of $\func$ of index $i$.

Further let $\bcomp$ be a connected component of $\partial\surf$.
Since $\func$ has no critical points on $\bcomp$, it follows that the gradient vector field $\grad\func$ is transversal to $\bcomp$ in every point $\pnt\in\bcomp$.

Consider the function $\sgn{\func}:\pi_0\partial\surf\to\{-1,1\}$ such that 
for every connected component $\bcomp$ of $\partial\surf$, regarded as an element of $\pi_0\partial\surf$, we have 
$\sgn{\func}(\bcomp)=+1$ provided 
$\grad\func$ is directed outward on all of $\bcomp$ and $\sgn{\func}(\bcomp)=-1$ otherwise.
We will call $\bcomp$ either \emph{$\func$-positive} or {\em $\func$-negative} in accordance with $\sgn{\func}(\bcomp)$.
Then the following quadruple:
\begin{equation}\label{equ:critical_type}
 \Crtype{\func} =
 \left(
    \crpntf{0}, \ \crpntf{1}, \ \crpntf{2}, \ \sgn{\func}
 \right)
\end{equation}
will be called the {\em critical type\/} of a Morse mapping $\func$.
Notice that the reversion of the orientation of $\onemanif$
interchanges $\crpntf{0}$ and $\crpntf{2}$ and replaces $\sgn{\bcomp}$ by $-\sgn{\bcomp}$.

Finally we will say that two Morse mappings $\func,\gfunc:\surf\to\onemanif$ are \emph{\sgmhom-homotopic} (belong to same connected path-component of $\MrsMP$), and write $\func\sgmh\gfunc$, if there is a continuous mapping $\Func:\surf\times\Interv\to\onemanif$ such that for every $t\in\Interv$ the function
$\func_t(x) = \Func(x,t) :\surf\to\onemanif$ is Morse.

\begin{theorem}[S.~Matveev~\cite{Kudr}, V.~Sharko~\cite{Sharko}]\label{th:MatvSharko}
Two Morse functions $\func,\gfunc:\surf\to\RRR$ are \sgmhom-homotopic iff $\Crtype{\func} = \Crtype{\gfunc}$.
Moreover, suppose that $\func=\gfunc$ is a neighborhood of some open-closed subset $V$ of $\partial\surf$. Then $\func\sgmh\gfunc$ with respect to some neighborhood of $V$.
\end{theorem}

The main result of this note is the following theorem:

\begin{theorem}\label{th:Main}
Two Morse functions $\func, \gfunc:\surf\to\Circle$ are \sgmhom-homotopic if and only if they are homotopic and $\Crtype{\func} = \Crtype{\gfunc}$.
\end{theorem}
The proof is heavily based on Theorem~\ref{th:MatvSharko} and the structure of minimal Morse functions on $\surf$.
\section{Preliminaries}
\label{sect:prelim}
We will regard $\Circle$ as $\RRR/\ZZZ$.
Let $\func:\surf\to\Circle$ be a Morse mapping.
Then a point $\pnt\in\Circle$ will be called a \emph{regular} value of $\func$ if $\func^{-1}(\pnt)$ contains no critical points of $\func$ and no connected components of $\partial\surf$.

The following construction will often be used.
Let $\apnt=0\in\Circle$ be a regular value of $\func$.
Let us cut $\surf$ along $\func^{-1}(0)$ and denote the obtained surface by $\csurf{\apnt}$. 
Let also $\cpr:\csurf{\apnt}\to\surf$ be the factor-map and
$\qpr(t) = e^{2\pi i t}:\RRR\to\Circle$ a universal covering.
Then there is a Morse function $\cfunc:\csurf{0}\to[0,1]$,
such that the following diagram is commutative:
\begin{equation}\label{equ:proj_tsurf}
\begin{CD}
\csurf{\apnt} @>{\cfunc}>> [0, 1] \\
@V{\cpr}VV  @V{\qpr}VV \\
\surf @>{\func}>> \Circle.
\end{CD}
\end{equation}

\subsection{Orientation of level-sets}
\label{subsect:orient_level_set}
Let $\func:\surf\to\Circle$ be a Morse mapping.
Since $\func$ is constant on components of $\partial\surf$ we have the following homomorphism 
$\func^{*}:H^{1}(\Circle)\to H^{1}(\surf,\partial\surf)$.
Let $\xi\in H^{1}(\Circle)\approx\ZZZ$ be the generator that yields chosen positive orientation of $\Circle$.
Then for every oriented closed curve $\omega:\Circle\to\surf$ we have 
\begin{equation}\label{equ:init_descr_int_ind}
   \deg(\func\circ\omega) = \func^{*}(\xi)(\omega).
\end{equation}

Suppose that $\surf$ is oriented.
Then there is an orientation of level-sets of $\func$ such that 
for every regular point $\apnt\in\surf$ of $\func$ and a tangent vector $v$ to $\func^{-1}\func(\apnt)$ at $\apnt$ the pair $(\grad\func(\apnt),v)$ gives the positive orientation of $T_{\apnt}M$.
Thus the level-sets of $\func$ can be regarded as elements of $H_1(\surf,\partial\surf)$.
 
Recall that there is an intersection form on $\surf$
$$
\langle\cdot,\cdot\rangle:
H_1(\surf,\partial\surf)\times H_1(\surf,\partial\surf) \to \ZZZ
$$
such that the correspondence $Z \mapsto \langle Z,\cdot\rangle$,  
$Z\in H_1(\surf,\partial\surf)$, yields an isomorphism $H_1(\surf,\partial\surf)\approx H^1(\surf,\partial\surf)$.
Then for $\pnt\in\Circle$ we have 
\begin{equation}\label{equ:descr_int_ind}
  \deg(\func\circ\omega) = \func^{*}(\xi)(\omega) = \langle \func^{-1}(\pnt),\omega\rangle. 
\end{equation}

\begin{lemma}\label{lm:homotopy_criterium}
Let $\func,\gfunc:\surf\to\Circle$ be two smooth functions which take constant values on connected components of $\surf$.
Then the following conditions are equivalent:
\begin{enumerate}
\item
$\func$ and $\gfunc$ are homotopic;
\item
$\func^{*} = \gfunc^{*}$;
\item
for every $\apnt,\bpnt\in\Circle$ the $1$-cycles $\func^{-1}(\apnt)$ and $\gfunc^{-1}(\bpnt)$ are homological in $H_{1}(\surf,\partial\surf)$.
\end{enumerate}
\end{lemma}
\proof
Equivalence (1)$\Leftrightarrow$(2) is well-known.

(2)$\Leftrightarrow$(3).
Let $\apnt,\bpnt\in\Circle$ and $\regcomp=\func^{-1}(\apnt)$ and  $\gregcomp=\gfunc^{-1}(\bpnt)$. 
Then $\regcomp=\gregcomp$ in $H_1(\surf,\partial\surf)$ iff 
$\langle\regcomp,\omega\rangle = \langle\gregcomp,\omega\rangle$ for every oriented closed curve in $\surf$.
In view of ~\eqref{equ:descr_int_ind} this is equivalent to the statement that $\func^{*}(\xi)=\gfunc^{*}(\xi)$.
\endproof

\subsection{Minimal Morse functions}
\label{subsect:minim_Morse_map}
Let $\compA$ and $\compB$ be two disjoint open-closed subsets of $\partial\surf$ (we do not require that $\compA\cup\compB=\partial\surf$).
Then $\compA$ and $\compB$ consist of connected components of $\partial\surf$. 

Recall that a Morse function $\func:\surf\to\RRR$ is \emph{minimal}, provided $\func$ has minimal number of critical points at each index among all Morse function on $\surf$.

The following statement is well-known, see e.g.~\cite{Wall:DiffTop}
\begin{lemma}\label{lm:minim_func}
Let $\sgn{}:\pi_0\partial\surf\to \{-1,1\}$ be an arbitrary function such that $\sgn{}(\compA)=0$ and $\sgn{}(\compB)=1$.
Then there exists a minimal Morse function $\func:\surf\to\imInterv$ 
such that $\func^{-1}(\aval)=\compA$, $\func^{-1}(\bval)=\compB$, and $\sgn{\func}=\sgn{}$.
Moreover, for every such a function we have

{\rm 1)} $\crpnt{\func}{0}=0$ provided $\sgn{}^{-1}(0)\not=\emptyset$; otherwise $\crpnt{\func}{0}=1$.

{\rm 2)} Similarly, if $\sgn{}^{-1}(1)\not=\emptyset$, then $\crpnt{\func}{2}=0$, otherwise $\crpnt{\func}{2}=1$.

Finally, every Morse function can be obtained from some minimal one by adding proper number of pairs of critical points of indexes $0$ and $1$ or $1$ and $2$.
\qed
\end{lemma}

\subsection{Unessential components}
\label{subsect:unessent_comp}
Let $\func:\surf\to\Circle$ be a Morse mapping, $\apnt$ a regular value of $\func$ and $\regcomp=\func^{-1}(\apnt)$.

A connected component $\concomp$ of $\surf\setminus\regcomp$ will be called \emph{essential} if either $\regcomp=\emptyset$ or $\func(\overline{\concomp})=\Circle$, otherwise $\concomp$ is \emph{unessential}.

Let $\concomp$ be an unessential component of $\surf\setminus\regcomp$.
Then $\concomp$ \emph{lower} if $\func(\overline{\concomp})\subset[\apnt, \apnt+d)$, for some $d \in (0,1)$.
Otherwise, $\func(\overline{\concomp})\subset(\apnt-d,\apnt]$ for some $d \in (0,1)$ and $\concomp$ will be called \emph{upper}.

Finally, we will say that $\func$ is \emph{$\apnt$-reduced}, provided all connected components of $\surf\setminus\regcomp$ are essential.

\begin{lemma}\label{prop:reduce_unnessent}
In the above notations, $\func$ is \sgmhom-homotopic to an $\apnt$-reduced Morse mapping.
\end{lemma}
\proof
Let $\concomp$ be an unessential component of $\surf\setminus\regcomp$.
We can assume that $\concomp$ is lower so that  $\func(\overline{\concomp})=[\apnt,d]$, where $0<\apnt<d<1$ and the interval $[0,\apnt]$ consists of regular values of $\func$ only.
Denote by $\aconcomp$ the connected component of $\func^{-1}[0,d]$ including $\concomp$.

Let also $\mu:[0,1]\to[0,1]$ be a $C^{\infty}$-function such that $\mu(0)=1$ and $\mu(1)=\apnt$.
Then it is easy to verify that the following mapping $\Func:\surf\times\Interv\to\Circle$ defined by
$$
\Func(\xpnt,t) = \left\{
\begin{array}{cl}
\mu(t) \func(\xpnt), & \xpnt\in\aconcomp\\
\func(\xpnt), & \xpnt\in\surf\setminus\concomp.
\end{array}
\right.
$$
is a \sgmhom-homotopy between $\func = \Func_0$ and the mapping $\gfunc=\Func_1$ such that
$\gfunc(\concomp) \subset [\apnt^2,\apnt d] \subset [0,\apnt)$, whence $\gfunc^{-1}(\apnt)= \regcomp \setminus \concomp$.

Then our lemma follows by the induction on the number of connected components of $\regcomp$.
\endproof

\subsection{Construction of Morse functions with given regular level-set}
\label{subsect:predict_level_set}
Let $\crv = \{\crv_1,\ldots,\crv_n\} \subset\Int\surf$ be a family of mutually disjoint two-sided simple closed curves, 
$\func:\surf\to\Circle$ a Morse function, $\apnt\in\Circle$ a regular value of $\func$, and $\regcomp=\func^{-1}(\apnt)$.

\begin{defn}\label{def:regular_family}
We will say that $\crv$ is \emph{$(\func,\apnt)$-regular} if 
$\regcomp\cap\crv=\emptyset$ and for every connected component $\concomp$ of $\surf\setminus(\regcomp\cup\crv)$ we have
\begin{enumerate}
\item[\rm(1)]
$\overline{\concomp}\cap\regcomp\not=\emptyset$ and $\overline{\concomp}\cap\crv\not=\emptyset$;
\item[\rm(2)]
$\func(\overline{\concomp})\not=\Circle$.
\end{enumerate}
\end{defn}

It follows from (2) that either $\func(\overline{\concomp})=[\apnt,\apnt+d]$ or $\func(\overline{\concomp})\subset[\apnt-d,\apnt]$ for some $d\in(0,1)$.
We will call $\concomp$ \emph{lower} in the first case and \emph{upper} in the second.
 
\begin{lemma}\label{lm:constr_hfunc}
Suppose that $\crv$ is $(\func,\apnt)$-regular.
Then there exists a Morse function $\hfunc:\surf\to\Circle$ such that 
\begin{enumerate}
\item
$\hfunc^{-1}(\apnt) = \regcomp = \func^{-1}(\apnt)$ and $\hfunc=\func$ near $\regcomp$;

\item
$\hfunc^{-1}(\bpnt)=\crv$ for some regular value $\bpnt$ of $\hfunc$;

\item 
$\Crtype{\hfunc} = \Crtype{\func}$.
\end{enumerate}
Then it follows from Lemma~\ref{lm:invf_invg} that 
$\hfunc \sgmh \func$ with respect to a neighborhood of $\regcomp$.
\end{lemma}
\proof
We can assume that $\apnt=0$ and $\bpnt=\frac{1}{2}$.
Let $\concomp$ be a connected component of $\surf\setminus(\regcomp\cup\crv)$.
If $\concomp$ is lower, then it follows from Definition~\ref{def:regular_family} and Lemma~\ref{lm:minim_func}, that there exists a \emph{minimal} Morse function 
$\hfunc_{\concomp}:\overline{\concomp}\to[0,\frac{1}{2}]$ such that 
$\hfunc_{\concomp}^{-1}(\frac{1}{2})=\overline{\concomp}\cap\crv$, 
$\hfunc_{\concomp}^{-1}(0)=\overline{\concomp}\cap\regcomp$, and $\sgn{\hfunc} = \sgn{\func}$.

Similarly, if $\concomp$ is upper, then we can construct a minimal Morse function  $\hfunc_{\concomp}:\overline{\concomp}\to[\frac{1}{2},1]$ such that $\hfunc_{\concomp}^{-1}(\frac{1}{2})=\overline{\concomp}\cap\crv$,
$\hfunc_{\concomp}^{-1}(1)=\overline{\concomp}\cap\regcomp$, and $\sgn{\hfunc} = \sgn{\func}$.


Then the union of all functions $\hfunc_{\concomp}$, where $\concomp$ runs all connected components of $\surf\setminus(\regcomp\cup\crv)$, gives 
a function $\hat\hfunc:\surf\to\Circle$  without critical points of indexes $0$ and $2$ and such that $\hat\hfunc^{-1}(0)=\regcomp$, $\hat\hfunc^{-1}(\frac{1}{2})=\crv$. 

Moreover, we can choose these functions so that $\hat\hfunc$ is smooth near $\regcomp\cup\crv$.
Then adding to $\tilde\hfunc$ a necessary number of pairs of critical points of indexes $0$ and $1$ or $1$ and $2$ we can obtain a Morse function $\hfunc:\surf\to\Circle$ satisfying the conditions (1)-(3) of our lemma.

Evidently, the condition (1) implies that $\hfunc$ is homotopic to $\func$, whence by Lemma~\ref{lm:invf_invg} we get $\func\sgmh\hfunc$ with respect to a neighborhood of $\regcomp$.
\endproof
\section{Proof of Theorem~\ref{th:Main}}
\label{sect:Proof_MainTheorem}
The necessity is obvious, therefore we will consider only sufficiency.
Let $\func, \gfunc:\surf\to\Circle$ be two Morse mappings that are homotopic and $\Crtype{\func}=\Crtype{\gfunc}$.
We have to show that $\func\sgmh\gfunc$.
First consider one particular case.
\begin{lemma}\label{lm:invf_invg}
Suppose that there exists a common regular value $\apnt$ of $\func$ and $\gfunc$ such that
$  \func^{-1}(\apnt) = \gfunc^{-1}(\apnt),$
and $\func = \gfunc$ in a neighborhood of $\func^{-1}(\apnt)$.
Then $\func\sgmh\gfunc$.
\end{lemma}
\proof
Denote $\regcomp = \func^{-1}(\apnt)$.
If $\regcomp=\emptyset$, then $\func$ and $\gfunc$ are mappings $\surf\to\Circle\setminus\{\apnt\}\approx \RRR$, whence by Theorem~\ref{th:MatvSharko}, $\func \sgmh \gfunc$.

Thus suppose that $\regcomp\not=\emptyset$ and let $\apnt=0$.
Then $\regcomp$ is a disjoint union of \emph{two-sided} simple closed curves.
Using the notations of~\eqref{equ:proj_tsurf},
we cut $\surf$ along $\regcomp$ and obtain liftings $\tfunc,\tgfunc:\csurf{\apnt}\to[0,1]$ of $\func$ and $\gfunc$ respectively.
Then
$\tfunc =\tgfunc$ near $\tregcomp = \pr^{-1}(\regcomp)$.

\begin{claim}
$\func$ is \sgmhom-homotopic to a Morse function $\hfunc:\surf\to\Circle$ such that $\Crtype{\hfunc|_{\overline{\concomp}}}=  \Crtype{\gfunc|_{\overline{\concomp}}}$
for every connected component $\concomp$ of $\surf\setminus\regcomp$.
\end{claim}
It follows from this claim that 
$\hfunc$ yields a Morse map $\thfunc:\tsurf\to[0,1]$ such that  
$\Crtype{\thfunc|_{\overline{\aconcomp}}}= \Crtype{\tgfunc|_{\overline{\aconcomp}}}$ for every connected component $\aconcomp$ of $\tsurf$.
Then from Theorem~\ref{th:MatvSharko} we obtain that  $\thfunc|_{\overline{\aconcomp}}\sgmh\tgfunc|_{\overline{\aconcomp}}$ with respect to a neighborhood of $\tregcomp\cap\overline{\aconcomp}$.
Hence $\thfunc\sgmh\tgfunc$ with respect to a neighborhood of $\tregcomp$ and 
therefore $\hfunc\sgmh\gfunc$ with respect to a neighborhood of $\regcomp$.
Thus $\func\sgmh\hfunc\sgmh\gfunc$. This will prove Lemma~\ref{lm:invf_invg}.

\proof[Proof of Claim]
It follows from Lemma~\ref{lm:minim_func} that such a function $\hfunc$ can be obtained from $\func$ by moving pairs of critical points of indexes $0$ and $1$ and indexes $1$ and $2$ from some connected components of $\surf\setminus\regcomp$ to another ones.

Consider the partition of $\surf$ by the connected components of level-sets of $\func$. Recall that the factor-space of $\surf$ by this partition admits a natural structure of a graph called \emph{Reeb} graph of $\func$.

\begin{center}
\begin{figure}[ht]
\includegraphics[width=0.8\textwidth]{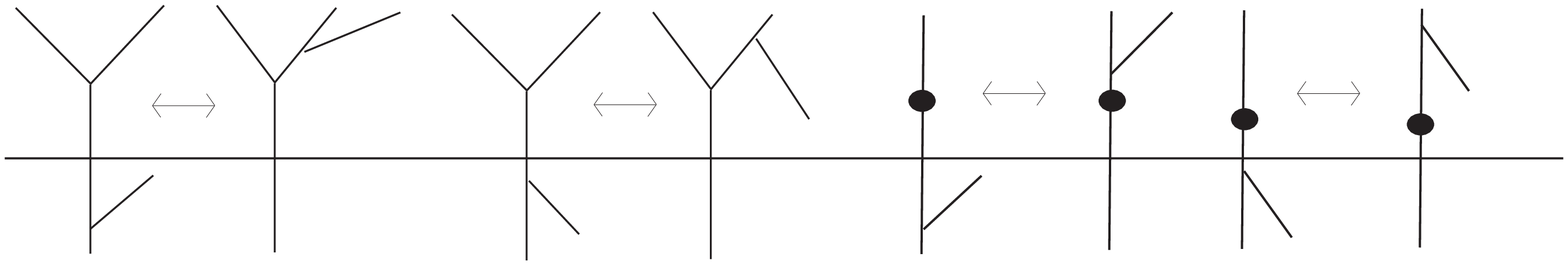}
\caption{}
\protect\label{fig:ReebGraph}
\end{figure}
\end{center}

Evidently, moves of pairs of critical points yield transformations of Reeb graph of $\func$ shifting edges with vertexes of degree $1$, see Figure~\ref{fig:ReebGraph}, where bold points are the vertices of degree $2$.

Every such a transformation can be realized by some \sgmhom-homotopy $\func_t, (t\in[0,1])$.

Moreover, let $\apnt$ be a value of of $\func$ corresponding to the level-set denoted in Figure~\ref{fig:ReebGraph} by long horizontal line.
Then $\func_t$ can be chosen so that $\func_0^{-1}(\apnt) = \func_1^{-1}(\apnt)$ and $\func_0=\func_1$ near $\func_0^{-1}(\apnt)$, while it is possible that 
$\func_0^{-1}(\apnt) \not= \func_t^{-1}(\apnt)$ for some $t\in(0,1)$.

Thus properly moving edges with vertexes of degree $1$ we can obtain from $\func$ a Morse function $\hfunc$ satisfying the statement of this lemma.
\endproof

Now Theorem~\ref{th:Main} is implied by the following two propositions and previous Lemma~\ref{lm:invf_invg}.

\begin{prop}\label{pr:XY_disj}
Let $\apnt$ and $\bpnt$ be regular values of $\func$ and $\gfunc$ respectively.
Suppose that 
$\func$ and $\gfunc$ are reduced with respect to $\apnt$ and $\bpnt$, and 
$\func^{-1}(\apnt) \cap \gfunc^{-1}(\bpnt)=\emptyset$.
Then $\func\sgmh\gfunc$.
\end{prop}

\begin{prop}\label{pr:DisjInv}
The functions $\func$ and $\gfunc$ are $\Sigma$-homotopic to Morse mappings 
$\func_{1}$ and $\gfunc_{1}$ respectively such that 
$\apnt$ and $\bpnt$ are regular values of $\func_1$ and $\gfunc_1$ respectively, and 
$\func_{1}^{-1}(\apnt) \cap \gfunc_{1}^{-1}(\bpnt) = \emptyset.$
Then by Proposition~\ref{pr:XY_disj}, $\func_1 \sgmh \gfunc_1$, whence $\func\sgmh\gfunc$.
\end{prop}

\section{Proof of Proposition~\ref{pr:XY_disj}}
\label{sect:Proof_Prop_Disj_to_Common}
Denote $\regcomp=\func^{-1}(\apnt)$ and $\gregcomp=\gfunc^{-1}(\bpnt)$, so we have $\regcomp\cap\gregcomp=\emptyset$.
It suffices to prove the following statement.

\begin{claim}\label{clm:Y_cuts_M_between_X0_X1}
Let $\concomp$ be a connected component of $\surf\setminus(\regcomp\cup\gregcomp)$.
Then 
$$\func(\overline{\concomp})\not=\Circle \qquad \text{and} \qquad \overline{\concomp}\cap\gregcomp\not=\emptyset.$$
Similarly, 
$\gfunc(\concomp)\not=\Circle$ and $\overline{\concomp}\cap\regcomp\not=\emptyset.$
Thus $\gregcomp$ is $(\func,\apnt)$-regular and $\regcomp$ is $(\gfunc,\bpnt)$-regular.
\end{claim}

It will follow from Lemma~\ref{lm:constr_hfunc} that there exists a Morse function $\hfunc$ such that $\hfunc^{-1}(\apnt)=\regcomp$, $\hfunc^{-1}(\bpnt)=\gregcomp$, 
$\hfunc=\func$ near $\regcomp$, and $\hfunc=\gfunc$ near $\gregcomp$.
Then by Lemma~\ref{lm:invf_invg} we will get $\func\sgmh\hfunc\sgmh\gfunc$.

\proof[Proof of Claim~\ref{clm:Y_cuts_M_between_X0_X1}]
(1) First suppose that $\surf$ is orientable.
Let us assume that $\apnt=0$.
Cutting $\surf$ along $\regcomp$ we obtain a connected surface $\csurf{0}$, the projection $\pr:\csurf{0}\to\surf$, and a Morse function $\tfunc:\csurf{0}\to[0,1]$ such that the commutative diagram~\eqref{equ:proj_tsurf} holds true.

Denote $\tregcomp_0 = \tfunc^{-1}(0)$, $\tregcomp_1 = \tfunc^{-1}(1)$, $\tgregcomp=\pr^{-1}(\gregcomp)$, and $\aconcomp=\pr^{-1}(\concomp)$.

Since $\func$ and $\gfunc$ are homotopic and $\surf$ is orientable, it follows from Lemma~\ref{lm:homotopy_criterium} that 
$1$-cycles $[\regcomp]$ and $[\gregcomp]$ are homological modulo $\partial\surf$. 
Whence $\tgregcomp$ separates $\tregcomp_0$ and $\tregcomp_1$ in $\csurf{0}$, i.e. for every connected subset $\nbh\subset\tsurf$ such that $\nbh\cap\tregcomp_i\not=\emptyset$ ($i=0,1$) we have that $\nbh\setminus\tgregcomp$ can be represented as a union of two disjoint open-closed subsets $\nbh_{i}$ such that $\nbh\cap\tregcomp_i\subset\nbh_{i}$. 

Suppose that $\func(\overline{\concomp})=\Circle$.
Then $\tfunc(\overline{\aconcomp})=[0,1]$.
Therefore $\overline{\aconcomp}\cap\tregcomp_i\not=\emptyset$ for $i=0,1$.
Hence $\overline{\aconcomp}\setminus\tgregcomp$, $\overline{\concomp}\setminus\gregcomp$, and therefore $\concomp\setminus\gregcomp$ are not connected which contradicts to the assumption.

If $\overline{\concomp}\cap\gregcomp=\emptyset$, then 
$\concomp$ is a connected component of $\surf\setminus\regcomp$.
But $\func(\overline{\concomp})\not=\Circle$, whence 
$\concomp$ is unessential for $(\func,\apnt)$, i.e. $\func$ is not $\apnt$-reduced.

(2) Suppose now that $\surf$ is non-orientable.
Let $\orpr:\orsurf\to\surf$ be the oriented double covering,
$\orfunc=\func\circ\orpr$, $\orgfunc=\gfunc\circ\orpr:\orsurf\to\Circle,$
$$
  \orregcomp=\orpr^{-1}(\regcomp)= \orfunc^{\,-1}(\apnt), \quad
  \orgregcomp=\orpr^{-1}(\gregcomp)= \orgfunc^{\,-1}(\bpnt), \quad \text{and} \quad \orconcomp=\orpr^{-1}(\concomp).
$$ 
If $\concomp$ is orientable, then $\orconcomp$ consists of two components 
 each homeomorphic to $\concomp$.
Otherwise, $\concomp$ is non-orientable, and $\orconcomp$ is an oriented double covering of $\concomp$.

Notice that $\orfunc$ is $\apnt$-reduced and $\orgfunc$ is $\bpnt$-reduced.
Indeed, for every connected component $\aconcomp$ of $\orsurf\setminus\orregcomp$ 
the set $\orpr(\aconcomp)$ is a connected component of $\surf\setminus\regcomp$, whence $\orfunc(\aconcomp) = \func( \orpr(\aconcomp) )\not=\Circle$.
The proof for $\orgfunc$ is similar.


Let $\oraconcomp$ be a connected component of $\orconcomp$.
Since $\orregcomp\cap\orgregcomp=\emptyset$, we get from orientable case of this claim that 
$$
\func(\concomp) = \orfunc(\oraconcomp)\not=\Circle \qquad  
\text{and} \qquad 
 \overline{\concomp}\cap\gregcomp \supset \orpr(\overline{\oraconcomp}\cap\orgregcomp)\not=\emptyset.
$$
This completes Claim~\ref{clm:Y_cuts_M_between_X0_X1} and Proposition~\ref{pr:XY_disj}.
\endproof
\section{Proof of Proposition~\ref{pr:DisjInv}}
\label{sect:Proof_Prop_Intersets_to_Disj}
We can assume that the intersection $\regcomp$ and $\gregcomp$ is transversal.
Let $n=\#(\regcomp\cap\gregcomp)$. If $n=0$, then our statement is just  Proposition~\ref{pr:XY_disj}.
Suppose that $n>0$.
We will show how to reduce the number of intersection points $\regcomp\cap\gregcomp$ by \sgmhom-homotopy.

For simplicity let $\apnt=0$.
We can also assume that $\func$ is $\apnt$-reduced.
Cutting $\surf$ along $\regcomp$ we obtain a connected surface $\csurf{0}$, the projection $\pr:\csurf{0}\to\surf$, and a Morse function $\tfunc:\csurf{0}\to[0,1]$ such that the commutative diagram~\eqref{equ:proj_tsurf} holds true.
Denote $\tregcomp_0 = \tfunc^{-1}(0)$, $\tregcomp_1 = \tfunc^{-1}(1)$, and $\tgregcomp=\pr^{-1}(\gregcomp)$.

Notice that $\tgregcomp$ consists of simple closed curves and arcs with ends at $\tregcomp$.
Let us divide $\tgregcomp$ by the following four disjoint subsets: $$\larcD, \ \larcU, \ \larcDU, \ \larcC,$$ where
$\larcD$ ($\larcU$) consists of arcs 
whose both ends belong to $\tregcomp_{0}$ ($\tregcomp_{1}$),
$\larcDU$ consists of arcs of $\tgregcomp$ connecting $\tregcomp_{0}$ with $\tregcomp_{1}$, and
$\larcC$ consists of simple closed curves of $\tgregcomp$.

Since $\func$ and $\gfunc$ are homotopic, it follows from~\eqref{equ:init_descr_int_ind}, that the restriction $\func|_{\gregcomp}$ is null-homotopic.
Hence $\larcD\not=\emptyset$ and $\larcU\not=\emptyset$.


Let $\nbh$ be a regular neighborhood of $ \tregcomp_0\cup\larcD$ such that $\partial\nbh$ does not intersect $\larcU\cup\larcC$ and transversely intersects every arc of $\larcDU$ at a unique point.

Let also $\aconcomps$ be the union of closures of those connected components $\aconcomp$ of $\tsurf\setminus(\tregcomp\cup\tcrv)$ for which 
$\overline{\aconcomp}\cap\tregcomp=\emptyset$.

Denote $\nbh'=\nbh\cup\aconcomps$, $\tcrv=\partial\nbh'$, and $\crv=\pr(\tcrv)$.
Then 
\begin{equation}\label{equ:intersect_decrease}
\#(\crv\cap\gregcomp) < \#(\regcomp\cap\gregcomp)=n.
\end{equation}
\begin{claim}\label{clm:crv_is_reg}
$\crv$ is $(\func,\apnt)$-regular. 
\end{claim}
Then by Lemma~\ref{lm:constr_hfunc} we can construct a Morse function $\hfunc:\surf\to\Circle$ such that $\hfunc^{-1}(0)=\regcomp$, $\hfunc^{-1}(\frac{1}{2})=\crv$, $\hfunc=\func$ near $\regcomp$, and $\hfunc\sgmh\func$.
It will follow from~\eqref{equ:intersect_decrease} by the induction on $n$ that $\hfunc\sgmh\gfunc$.

\proof[Proof of Claim~\ref{clm:crv_is_reg}]
We have to show that for every connected component $\concomp$ of $\surf\setminus(\regcomp\cup\crv)$ the following conditions hold true:
\begin{equation}\label{equ:cond_for_concomp}
\func(\overline{\concomp})\not=\Circle, \qquad \overline{\concomp}\cap\crv\not=\emptyset, \qquad \overline{\concomp}\cap\regcomp\not=\emptyset.
\end{equation}
Denote $\aconcomp=\pr^{-1}(\concomp)$.
Then~\eqref{equ:cond_for_concomp} 
are equivalent to the following ones:
$$
\tfunc(\overline{\aconcomp})\not=[0,1], \quad 
\overline{\aconcomp}\cap\tcrv\not=\emptyset, \quad \overline{\aconcomp}\cap\tregcomp\not=\emptyset.
$$

Since $\partial\nbh$ separates $\tregcomp_0$ and $\tregcomp_1$ in $\tsurf$, we get $\tfunc(\overline{\aconcomp})\not=[0,1]$.

Moreover, if $\overline{\aconcomp}\cap\tcrv=\emptyset$, then $\aconcomp$ is in fact a connected component of $\tsurf$. But $\tfunc(\overline{\aconcomp})\not=[0,1]$ implies $\func(\overline{\concomp})\not=\Circle$, whence $\concomp$ is unessential with respect to $(\func,\apnt)$ which contradicts to the assumption.

Suppose that $\overline{\aconcomp}\cap\tregcomp=\emptyset$. Then $\aconcomp$ is a connected component of $\tsurf\setminus\crv$. 
Therefore, $\overline{\aconcomp}\cap\nbh\subset\partial\nbh$, whence $\aconcomp\subset\aconcomps$. But every connected component of $\aconcomps$ evidently intersect $\tregcomp_0$, thus $\aconcomp\cap\tregcomp_0\not=\emptyset$ which contradicts to the assumption.
\endproof

\section{Acknoledgements}
I am sincerely grateful to V.~V.~Sharko for help and attention to my work. I also want to thank my colleagues M.~A.~Pankov, E.~A.~Polulyah, A.~O.~Prishlyak, and I.~Y.~Vlasenko for many useful discussions.

\end{document}